\documentclass[10pt]{amsart}
\input{epsf}
\usepackage{amssymb,latexsym}
\topskip  \newcommand{\emptyword}{\varepsilon}
\newcommand{\resp}[1]{{\rm [}respectively; #1{\rm ]}}

\newcommand{\mxgddyckI}{1234,\linebreak[0]1243,\linebreak[0]1423,\linebreak[0]4123}
\newcommand{\mxgddyckII}{1324,\linebreak[0]1342,\linebreak[0]1432,\linebreak[0]4132}
\newcommand{\mxgddyckIII}{2134,\linebreak[0]2143,\linebreak[0]2413,\linebreak[0]4213}

\newcommand{\mxgddyckV}{2314,\linebreak[0]2413,\linebreak[0]3142,\linebreak[0]3241}

\newcommand{\mxgddyckVII}{1234,\linebreak[0]1324,\linebreak[0]2134,\linebreak[0]2314}
\newcommand{\mxgddyckVIII}{1234,\linebreak[0]2134,\linebreak[0]2314,\linebreak[0]3124}
\newcommand{\mxgddyckIX}{1324,\linebreak[0]2134,\linebreak[0]2314,\linebreak[0]3124}
\newcommand{\mxgddyckX}{1324,\linebreak[0]2134,\linebreak[0]3124,\linebreak[0]3214}
\newcommand{\mxgddyckXI}{1324,\linebreak[0]2314,\linebreak[0]3124,\linebreak[0]3214}
\newcommand{\mxgddyckXII}{1342,\linebreak[0]2341,\linebreak[0]3142,\linebreak[0]3241}

\newcommand{\mxgddyckXIV}{1324,\linebreak[0]1342,\linebreak[0]2314,\linebreak[0]2341}
\newcommand{\mxgddyckXV}{1342,\linebreak[0]2341,\linebreak[0]2431,\linebreak[0]3241}

\numberwithin{equation}{section}
\newtheorem{theorem}{Theorem}[section]
\newtheorem{proposition}[theorem]{Proposition}

\newtheorem{definition}[theorem]{Definition}

\newtheorem{remark}[theorem]{Remark}
\newtheorem{lemma}[theorem]{Lemma}

\newtheorem{example}[theorem]{Example}

\begin{document}
\def\mn{\mbox{-}}
\def\re{R^e}
\def\ro{R^o}
\def\ie{I^e}
\def\io{I^o}
\def\je{J^e}
\def\jo{J^o}
\def\SS{\mathfrak S}
\def\II{\mathfrak I}
\def\JJ{\mathfrak L}
\def\ax{\overline{x}}
\def\ttx{\left(\frac{1}{2\sqrt{x}}\right)}
\def\txx{\left(\frac{1}{2x}\right)}

\title{Some statistics on restricted $132$ involutions}

\author[Olivier Guibert and Toufik Mansour]
{O. Guibert and T. Mansour}


\maketitle

\begin{center}
{\small LaBRI (UMR 5800), Universit\'e Bordeaux 1,
       351 cours de la Lib\'eration,
       33405 Talence Cedex, France}

{\tt guibert@labri.fr}, {\tt toufik@labri.fr}
\end{center}

\section*{Abstract}
In \cite{GM} Guibert and Mansour studied involutions on $n$
letters avoiding (or containing exactly once) $132$ and avoiding
(or containing exactly once) an arbitrary pattern on $k$ letters.
They also established a bijection between $132$-avoiding
involutions and Dyck word prefixes of same length. Extending this
bijection to bilateral words allows to determine more parameters;
in particular, we consider the number of inversions and rises of
the involutions onto the words. This is the starting point for
considering two different directions: even/odd involutions and
statistics of some generalized patterns. Thus we first study
generating functions for the number of even or odd involutions on
$n$ letters avoiding (or containing exactly once) $132$ and
avoiding (or containing exactly once) an arbitrary pattern $\tau$
on $k$ letters. In several interesting cases the generating
function depends only on $k$ and is expressed via Chebyshev
polynomials of the second kind. Next, we consider other statistics
on $132$-avoiding involutions by counting an occurrences of some
generalized patterns, related to the enumeration according to the
number of rises.

\section{Introduction}
The aim of this paper is to give analogies of enumerative results
on certain classes of permutations characterized by
pattern-avoidance in the symmetric group $\SS_n$. In
$\II_n=\{\pi\in \SS_n : \pi=\pi^{-1}\}$ we identify classes of
restricted involutions with enumerative properties analogous to
results for permutations. More precisely, we study generating functions
for the number of even or odd involutions in $\II_n$ avoiding (or
containing exactly once) $132$, and avoiding (or containing
exactly once) an arbitrary permutation $\tau\in\SS_k$. Moreover we
consider statistics of some generalized patterns. In the remainder
of this section we present a brief account of earlier works which
motivated our investigation, basic definitions used throughout
the paper, and the organization of this paper.

\subsection{Background}
\label{intro-back}
In this subsection we present some classical statistics on permutations,
pattern avoidance for permutations,
some earlier results, generalized patterns, and pattern avoidance
for involutions.

Let $\pi\in\SS_n$. The number of {\em inversions} of $\pi$ is
given by $|\{(i,j): \pi_i>\pi_j,\ 1\leq i<j\leq n\}|$.
We say $\pi$ is an {\em even permutation} \resp{{\em
odd permutation}} if $\pi$ is a permutation together with even
\resp{odd} number of inversions.
We say two consecutive elements
$\pi_j$ and $\pi_{j+1}$ form a {\em rise} \resp{{\em descent}} if
$\pi_j<\pi_{j+1}$ \resp{$\pi_j>\pi_{j+1}$} where $1\leq j\leq
n-1$.
We say $\pi_j$ is a {\em right-to-left maximum} \resp{{\em
left-to-right minimum}} if $\pi_j>\pi_i$ \resp{$\pi_j<\pi_i$} for
all $j<i\leq n$ \resp{$1\leq i<j$} where $1\leq j\leq n$.
The number of {\em fixed points} of $\pi$ is given by
$|\{j: \pi_j=j,\ 1\leq j\leq n\}|$.

Let $\pi\in\SS_n$ and $\tau\in\SS_k$ be two permutations. We say
that $\pi$ {\em contains\/} $\tau$ if there exists a subsequence
$1\leq i_1<i_2<\dots<i_k\leq n$ such that
$(\pi_{i_1},\dots,\pi_{i_k})$ is order-isomorphic to $\tau$; in
such a context $\tau$ is usually called a {\em pattern\/}. We say
that $\pi$ {\em avoids\/} $\tau$, or is $\tau$-{\em avoiding\/},
if such a subsequence does not exist. The set of all
$\tau$-avoiding permutations in $\SS_n$ is denoted $\SS_n(\tau)$.
For an arbitrary finite collection of patterns $T$, we say that
$\pi$ avoids $T$ if $\pi$ avoids any $\tau\in T$; the
corresponding subset of $\SS_n$ is denoted $\SS_n(T)$. For
example, $34521\in\SS_5(132)$ whereas $31542\not\in\SS_5(132)$
because it contains four subsequences (that is $354$, $154$,
$152$, and $142$) order-isomorphic to $132$.

While the case of permutations avoiding a single pattern has
attracted much attention, the case of multiple pattern avoidance
remains less investigated. In particular, it is natural, as the
next step, to consider permutations avoiding pairs of patterns
$\tau_1$, $\tau_2$. This problem was solved completely for
$\tau_1,\tau_2\in\SS_3$ (see \cite{SS}), for $\tau_1\in\SS_3$ and
$\tau_2\in\SS_4$ (see \cite{W}), and for $\tau_1,\tau_2\in\SS_4$
(see \cite{B1,Km} and references therein). Several recent papers
\cite{CW,MV1,Kr,MV2,MV3,MV5} deal with the case $\tau_1\in\SS_3$,
$\tau_2\in\SS_k$ for various pairs $\tau_1,\tau_2$. Another
natural question is to study permutations avoiding $\tau_1$ and
containing $\tau_2$ exactly $t$ times. Such a problem for certain
$\tau_1,\tau_2\in\SS_3$ and $t=1$ was investigated in \cite{Ro},
and for certain $\tau_1\in\SS_3$, $\tau_2\in\SS_k$ in
\cite{RWZ,MV1,Kr,MV3,MV5}.

In \cite{BS} Babson and Steingr\`{\i}msson introduced {\em
generalized patterns} that allow the requirement that two adjacent
letters in a pattern must be adjacent in the permutation. In this
context, we write a classical pattern with dashes between any two
adjacent letters of the pattern, say $1243$, as $1\mn2\mn4\mn3$,
and if we write, say $12\mn4\mn3$, then we mean that if this
pattern occurs in permutation $\pi$, then the letters in the
permutation $\pi$ that correspond to $1$ and $2$ are adjacent. For
example, the permutation $\pi=3542176$ has only one occurrence of
the pattern $12\mn43$, namely the subsequences $3576$, whereas
$\pi$ has two occurrences of the pattern $1\mn2\mn4\mn3$, namely
the subsequences $3576$, and $3476$. Claesson \cite{C} presented a
complete solution for the number of permutations avoiding any
single $3$-letters generalized pattern with exactly one
adjacent pair of letters. Claesson and Mansour \cite{CM} (see also
\cite{M1,M2,MS}) presented a complete solution for the number of
permutations avoiding any double $3$-letters generalized patterns
with exactly one adjacent pair of letters. Besides, Kitaev
\cite{Ki} investigated simultaneous avoidance of two or more
$3$-letters generalized patterns without internal dashes.

An {\em involution} $\pi$ is a permutation such that
$\pi=\pi^{-1}$ (permutation equivalent to its inverse). Some authors
considered involutions avoiding patterns. In particular, someones
studied the enumeration of involutions in $\II_n(12\ldots k)$ in general
\cite{Regev,Gessel} or for some $k$ \cite{GouyouBYoung} because
this pattern is directly connected to Young tableaux of bounded
height by Robinson-Schensted algorithm \cite{Robinson,Schensted}.
Moreover, some other ones
\cite{GuibertThese,GPP} considered several sets of avoiding
involutions enumerated by the Motzkin numbers.  It already remains
a connected open problem: in \cite{GuibertThese} conjectures that
$\II_n(1432)$ is also enumerated by the $n$th Motzkin number.

\subsection{Basic tools}
\label{intro-tools} In this subsection we present some tools we
will use later as Chebyshev polynomials, Dyck word prefixes and
bilateral words, and generating trees. As an example, we also
present a sorting algorithm.

{\em Chebyshev polynomials of the second kind\/} (in what follows
just Chebyshev polynomials) are defined by
\begin{equation}
U_r(\cos\theta)=\frac{\sin(r+1)\theta}{\sin\theta}
\end{equation}
for $r\geq0$. The Chebyshev polynomials satisfy the following
recurrence $U_r(t)=2tU_{r-1}(t)-U_{r-2}(t)$ for $r\geq2$ together
with $U_0(t)=1$ and $U_1(t)=2t$. Evidently, $U_r(x)$ is a
polynomial of degree $r$ in $x$ with integer coefficients.
Chebyshev polynomials were invented for the needs of approximation
theory, but are also widely used in various other branches of
mathematics, including algebra, combinatorics, and number theory
(see \cite{Ri}). Apparently, for the first time the relation
between restricted permutations and Chebyshev polynomials was
discovered  by Chow and West in \cite{CW}, and later by Mansour
and Vainshtein \cite{MV1,MV2,MV3,MV5} and Krattenthaler \cite{Kr}.
These results are related to a rational function
\begin{equation}
R_k(x)=\frac{U_{k-1}\ttx}{\sqrt{x}U_k\ttx} \label{rk}
\end{equation}
for all $k\geq 1$. For example, $R_1(x)=1$,
$R_2(x)=\frac{1}{1-x}$, and $R_3(x)=\frac{1-x}{1-2x}$. It is easy
to see that for any $k$, $R_k(x)$ is rational in $x$ and satisfies
the following equation (see \cite{MV1,MV3,MV5})
\begin{equation}
R_k(x)=\frac{1}{1-xR_{k-1}(x)}. \label{eqrk}
\end{equation}
For all $k\geq 0$, we define
\begin{equation}
\re_k(x)=\frac{1}{2}(R_{k}(x)+R_k(-x)) \qquad\mbox{and}\qquad
\ro_k(x)=\frac{1}{2}(R_k(x)-R_k(-x)). \label{eqeo}
\end{equation}

Let $L$ be a set of letters; the set of all words on $L$ we denote
by $L^*$. The length of a word $w$ we denote by $|w|$ is its
number of letters. The number of occurrences of the letter $l\in
L$ in a word $w$ we denote by $|w|_l=|\{w_i : w_i=l\}|$. We denote
by $\emptyword$ the empty word (of length $0$).
\\
{\em Dyck word prefixes} are words $w$ in $\{x,\ax\}^*$ such that
for all $w = w' w''$, $|w'|_x \geq |w'|_{\ax}$. The set of all
Dyck word prefixes we denote by $P_{x,\ax}$. For example, $xxxx$,
$xxx\ax$, $xx\ax x$, $xx\ax\ax$, $x\ax xx$, and $x\ax x\ax$ are
all the Dyck word prefixes of length $4$. Dyck word prefixes in
$P_{x,\ax}$ of length $n$ are enumerated by the central binomial
coefficient $a_n=\binom{n}{[n/2]}$ for all $n \geq 0$.
A Dyck word prefix corresponds to the beginning of a Dyck word.
\\
{\em Dyck words} are Dyck word prefixes $w$ in $\{x,\ax\}^*$ such
that $|w|_x=|w|_{\ax}$. For example, $xx\ax\ax$ and $x\ax x\ax$
are all the Dyck words of length $4$. Dyck words of length $2n$ are
enumerated by $C_n=\frac{1}{n+1}\binom{2n}{n}$ the $n$th Catalan
number whose generating function is $C(x)=\frac{1-\sqrt{1-4x}}{2x}$.
\\
{\em Bilateral words} are words in $\{x,\ax\}^*$ such that $|w|_x
= |w|_{\ax}$ or $|w|_x = |w|_{\ax} - 1$.
The set of all bilateral words we denote by $B_{x,\ax}$.
For example, $\ax\ax xx$, $\ax xx\ax$, $\ax x\ax x$, $xx\ax\ax$,
$x\ax\ax x$, and $x\ax x\ax$ are all the bilateral words of length $4$.

We note $\Xi$ the well known bijection which is an example of a
result due to Chottin and Cori~\cite{CC} (see also \cite{GM}, just
before Theorem~2.4) between a Dyck word prefix $w_0 x w_1 x \ldots
x w_p\in P_{x,\overline{x}}$ where $w_i$ is a Dyck word for all $0
\leq i \leq p$ and the bilateral word $w_0 \overline{x} w_1 \ax
\ldots \ax w_{[(p-1)/2]} \ax w_{[(p+1)/2]} x w_{[(p+3)/2]} x
\ldots x w_p \linebreak[0] \in B_{x,\ax}$ of the same length. For example, the
Dyck word prefix $x\ax xxx\ax x\ax xxx\ax\ax x$ is in bijection by
$\Xi$ with the bilateral word $x\ax\ax\ax x\ax x\ax xxx\ax\ax x$.
Moreover, the words of length $4$ given in the previous paragraph
are respectively in bijection. So, the number of Dyck word
prefixes of length $n$ is equal to the number of bilateral words of
length $n$ trivially enumerated by the central binomial
coefficient $a_n$.

Following \cite{Chung} (see also several PHD-thesis
\cite{GireThese,GuibertThese,PergolaThese,WestPHD}), a {\em
generating tree} of a set of objects is a tree subject to the
conditions that each object of length $n$ appears once and only
once on a vextex of level $n$ and that the edges correspond to a
manner to grow the objects. In order to characterize a generating
tree by a {\em succession system} we associate to each object a
label such that any two nodes have the same label if their
subtrees are isomorphic. Therefore it suffices to specify the
label of the root and a set of {\em succession rules} explaining
how to derive from the label of a parent the labels of all of its
children.
For example, in \cite[Theorem~2.4]{GM}, Guibert and Mansour
established a bijection called $\Phi$ between $\II_n(132)$
and Dyck word prefixes in $P_{x,\ax}$ of length $n$.
These objects can both be characterized by a succession system
whose root is $(0)$, and whose succession rules are $(0)\leadsto (1)$
and $(p)\leadsto (p+1),(p-1)$ if $p\geq1$. They also stated that
the number of fixed points of the involution corresponds to the
difference between the number of letters $x$ and $\overline{x}$
into the Dyck word prefix. So, they established
\cite[Corollary~2.5]{GM} that the number of $132$-avoiding
involutions of length $n$ having $p$ fixed points with $0 \leq p
\leq n$ (and $p$ is odd if and only if $n$ is odd) is the ballot
(or Delannoy) number $a_{n,p} = {n \choose \frac{n+p}{2}} - {n
\choose \frac{n+p}{2}+1}$.

One very nice proof that $|\SS_n(132)|=C_n$ and that the numbers
of  $132$-avoiding permutations of length $n$ having $s$
right-to-left maxima is equal to the ballot number
$a_{2n-s-1,s-1}={2n-s-1 \choose n-1}-{2n-s-1 \choose n}$ is due to
Knuth~\cite{Kn}. It consists in using a sorting algorithm (the
output is the mirror of the identity permutation) with one stack
(verifying its elements decrease from the top to the basis). This
defines a bijection with these one stack-sortable permutations (by
taking successively the last element, the previous one, until the
first one) of length $n$ having $s$ right-to-left maxima and Dyck
words (coding the movements of insertion/deletion on the stack) of
length $2n$ having $s$ primitive Dyck words (that is the Dyck
paths enumerated according to the number of times they touch the
horizontal axis). For example, to $43512 \in \SS_5(132)$ having
two right-to-left maxima, $5$ and $2$, corresponds the Dyck word
$xx\ax\ax xx\ax x\ax\ax$ (insert $2$, insert $1$, delete $1$,
delete $2$, insert $5$, insert $3$, delete $3$, insert $4$, and
then delete $4$ and $5$ to empty the stack), concatening the two
primitive Dyck words $xx\ax\ax$ and $xx\ax x\ax\ax$.

We also encounter several times the $n$th Fibonacci number
$F_n$ with $F_0=F_1=1$ and $F_{n}=F_{n-1}+F_{n-2}$ whose
generating function is $F(x)=\frac{1}{1-x-x^2}$.

\subsection{Organization of the paper}
In Section~\ref{sec2} we establish that the correspondence $\Phi
\circ \Xi$ between $132$-avoiding involutions and bilateral words
allows to determine more parameters. In particular, we consider
the number of inversions (and also the number of even or odd
involutions) and the number of rises of the involutions onto the
words.
\\
In Section~\ref{sec3} we consider the four cases of even or odd
involutions avoiding (or containing exactly once) $132$ and
avoiding (or containing exactly once) an arbitrary pattern $\tau$
on $k$ letters. In several interesting cases the generating
function depend only on $k$ and is expressed via Chebyshev
polynomials of the second kind.
\\
Finally, in Section~\ref{sec4} we present other statistics on
$132$-avoiding involutions by considering the distribution
according to the number of occurrences of some generalized
patterns. In particular, we relate some of these results to
Subsection~\ref{sec2Phi} because the number of rises of a
permutation is, evidently, given by the number of occurrences of
the generalized pattern $12$.
\section{Restricted $132$ involutions: number of inversions, and number of rises}
\label{sec2}

This section presents refinements of bijections $\Phi$ and $\Psi$
given in \cite{GM}. This allows us to determine the number of even
\resp{odd} involutions in $\II_n$ avoiding (or containing exactly
once) $132$. Besides, we determine the number of
involutions in $\II_n(132)$ having exactly $r$ rises.
\subsection{Refinements of bijection $\Phi$}
\label{sec2Phi}

In this subsection, following to \cite{GM}, we establish the
correspondence $\Phi \circ \Xi$ between $132$-avoiding involutions
and bilateral words which allows to determine more parameters. In
particular, we consider the number of inversions and of rises of the
involutions onto the words. So we recall the bijection $\Phi$
given in \cite[Theorem~2.4]{GM}.

Let $\pi\in\II_n(132)$ having $p$ fixed points. We have $\pi =
\pi' \pi'' x \pi'''$ with $|\pi'|=\frac{n-p}{2}$ ($\pi'$ has no
fixed points and constitutes cycles with $\pi''$ or $\pi'''$),
$\pi''$ does not contain fixed point ($\pi''$ constitutes cycles
with $\pi'$) and $\pi(x)=x$ ($x$ is the first fixed point and
$\pi'''$ constitutes cycles with $\pi'$ and/or contains fixed
points). We obtain two involutions in $\II_{n+1}(132)$ from $\pi$:
the first one is given by inserting a fixed point between $\pi'$
and $\pi''$, and the second one (if and only if $\pi$ has at least
one fixed point) is given by modifying the first fixed point $x$
by a cycle starting between $\pi'$ and $\pi''$.

Let $w\in P_{x,\ax}$ of length $n$ such that $|w|_x-|w|_{\ax}=p$.
So we have $w = w_0 x w_1 x \ldots x w_p$ where $w_i$ is a Dyck
word for all $0 \leq i \leq p$. We obtain two Dyck word prefixes
of length $n+1$ from $w$: $x w$ and $x w_0 \ax w_1 x \ldots x w_p$
(the latter if $p>0$).

Of course, the same construction can be applied to bilateral words
in bijection by $\Xi$ with Dyck word prefixes. We obtain two
bilateral words of length $n+1$ from $w_0 \ax w_1 \ldots \ax
w_{[(p+1)/2]} x \ldots x w_p$ in $B_{x,\ax}$ of length $n$: $\ax
w_0 \ax w_1 \ldots \ax w_{[p/2]} x \ldots x w_p$ and $x w_0 \ax
w_1 \ldots \ax w_{[p/2]+1} x \ldots x w_p$ (the latter if $p>0$).

The three generating trees for $132$-avoiding involutions, Dyck
word prefixes, and bilateral words are characterized by the same
succession system given in \cite[between Theorem~2.4 and
Corollary~2.5]{GM}:
\begin{equation}
\left\{\begin{array}{lcll}
    \multicolumn{4}{l}{(0)} \\
    (0) & \leadsto & (1) & \\
    (p) & \leadsto & (p+1) , (p-1) & \mbox{if } p \geq 1
    \end{array}\right..
\label{S}
\end{equation}
\begin{definition}
\label{def-invers} Let $w = w_0 x w_1 x w_2 x \ldots x w_p$ a Dyck
word prefix where $w_i$ is a Dyck word for all $0 \leq i \leq p$.
We define $i(w)$ by $i(\emptyword)=0$ and $$i(w) =
\left\{\begin{array}{lll}
    \multicolumn{2}{l}{|w_1 w_2 \ldots w_p| + i(w_1 x w_2 x \ldots x w_p)} \\
    &   \mbox{if } w_0 = \emptyword \\
    \multicolumn{2}{l}{
    |w_0' \ax w_0'' w_1 w_2 \ldots w_p| +
    |w_0'' w_1 w_2 \ldots w_p| +
    i(w_0' w_0'' x w_1 x w_2 x \ldots x w_p)
    } \\
    &   \mbox{if } w_0 = x w_0' \ax w_0''
        \mbox{ where } w_0' \mbox{ and } w_0'' \mbox{ are Dyck words}.
    \end{array}\right..$$
In such a context $i(w)$ can be called the number of {\em right
Dyck steps} of $w$.
\end{definition}

Of course, this statistic on Dyck word prefixes is similarly
defined on bilateral words.

\begin{theorem}
Let $\pi$ be an $132$-avoiding involution and $w$ be the bilateral
word in correspondence by $\Phi \circ \Xi$. Moreover, if $\pi$ has
length $n$, $p$ fixed points, $i$ inversions, and $r$ rises, then
$w$ has length $n$, $[(p+1)/2]$ minimal nonpositive height,
$i=i(w)$ right Dyck steps, and
$r=|wx|_{xx}+|w|_{\overline{x}\overline{x}}$ double identical
consecutive steps.
\end{theorem}
\begin{proof}
In \cite{GM} proved that Succession system~\ref{S} established a
one-to-one corespondence between involutions $\pi \in \II_n(132)$
having $p$ fixed points and Dyck word prefixes $w = w_0 x w_1 x
\ldots x w_p$ of length $n$ where $w_i$ is a Dyck word for all $0
\leq i \leq p$. So it is sufficient to consider the two others
parameters, inversions and rises. We show that Succession
system~\ref{S} can be extended to include these two parameters. We
consider the decomposition $\pi = \pi' \pi'' x \pi'''$ ($x$ is the
first fixed point and $\pi''$ contains only cycles connected to
$\pi'$).

First, we consider the number of inversions $i$ of involutions
that is the number of right Dyck steps of the Dyck word prefixes
(or equivalently bilateral words) and we obtain the succession
system:
\begin{equation}
\left\{\begin{array}{lcll}
    \multicolumn{4}{l}{(0,0)} \\
    (0,i) & \leadsto & (1,i+n-p) & \\
    (p,i) & \leadsto & (p+1,i+n-p) , (p-1,i+n-p+1) & \mbox{if } p \geq 1
    \end{array}\right..
\label{SS}
\end{equation}
Of course, if $n=0$, we have $i=0$ for the empty involution and
word.
\\
The involution obtained by inserting a fixed point between $\pi'$
and $\pi''$ has $i$ inversions unchanged from $\pi$ and
$\frac{n-p}{2}$ more for the new fixed point and $\frac{n-p}{2}$
more for $\pi'$ whereas the involution obtained by transforming
the first fixed point $x$ by a cycle starting between $\pi'$ and
$\pi''$ has $i$ inversions unchanged from $\pi$ and
$\frac{n-p}{2}$ more for $\pi'$ and $1+\frac{n-p}{2}$ more for the
new element between $\pi'$ and $\pi''$.

We now consider the Dyck word prefix (or equivalently the
bilateral word). By Definition~\ref{def-invers} we have
$$i(xw)=|w_0 w_1\ldots w_p|+i(w)=n-p+i(w)$$
for the case where $p$ increases and
$$\begin{array}{l}
i(x w_0 \overline{x} w_1 x w_2 x \ldots x w_p) = |w_0\overline{x}
w_1 w_2 \ldots w_p| + |w_1 w_2 \ldots w_p| + i(w_0 w_1 x w_2 x
\ldots x w_p)\\
\qquad\qquad\qquad\qquad\qquad\,\,= 1 + n - p + |w_1 w_2 \ldots
w_p| + i(w) - |w_1 w_2 \ldots w_p|\\
\qquad\qquad\qquad\qquad\qquad\,\,= 1+n-p+i(w) \end{array}$$ for
the case where $p$ decreases.

Secondly, we consider the number of rises $r$ of involutions that
is the number of double identical consecutive steps of bilateral
words. We must add an extra technical parameter $b$ (an indicator
$0$ or $1$) and we obtain the succession system:
\begin{equation}
\left\{\begin{array}{lcll}
    \multicolumn{4}{l}{(0,0,0)} \\
    (0,r,b) & \leadsto & (1,r+b,1) & \\
    (p,r,b) & \leadsto & (p+1,r+b,1) , (p-1,r+1-b,0) & \mbox{if } p \geq 1
    \end{array}\right.
\label{SSS}
\end{equation}
with $b=1$ \resp{$0$} if $\pi''=\emptyword$ \resp{$\neq
\emptyword$} and $b=1$ \resp{$0$} if $w_0=\emptyword$ \resp{$\neq
\emptyword$} that is $w$ starts by $\overline{x}$ \resp{$x$}.
\\
Of course, if $n=0$, we have $r=0$ and $b=0$ for the empty
involution and word.
\\
If $\pi''$ is \resp{is not} empty and $\pi$ has $r$ rises then the
involution obtained by inserting a fixed point has $r+1$
\resp{$r$} rises (one rise \resp{descent} more just after the new
fixed point) and $b$ is still equal to $1$ whereas the involution
obtained by transforming the first fixed point by a cycle has $r$
\resp{$r+1$} rises (one descent \resp{rise} more just after
\resp{before} the new element belonging to the new cycle) and $b$
becomes $0$.

We now consider the bilateral word. If $w_0$ is \resp{is not}
empty and $w$ is such that
$|wx|_{xx}+|w|_{\overline{x}\overline{x}}=r$, then the new
bilateral word has $r+1$ \resp{$r$} double identical consecutive
steps (it starts by $\overline{x}\overline{x}$
\resp{$\overline{x}w_0=\overline{x}xw_0'}$) and $b$ is still equal
to $1$ for the case where $p$ increases whereas the new bilateral
word has $r$ \resp{$r+1$} double identical consecutive steps (it
starts by $x\overline{x}w_1$ \resp{$xw_0\overline{x}w_1$}) and $b$
becomes $0$ for the case where $p$ decreases.
\end{proof}

Let $a'(n,p,i)$ be the number of involutions in $\II_n(132)$
having $p$ fixed points and $i$ inversions. Immediately, we deduce
from Succession system~\ref{SS} that
    $$a'(n,p,i)=a'(n-1,p-1,i+p-n)+a'(n-1,p+1,i+p-n+1)$$
together with $a'(0,0,0)=1$.
\begin{theorem}
\label{even132} The number of even \resp{odd} involutions in
$\II_n(132)$ is given by
\begin{center}
$\binom{n-1}{2[(n+1)/4]}$ $\left[ {respectively;}\
\binom{n-1}{1+2[(n-2)/4]} \right].$
\end{center}
\end{theorem}
\begin{proof}
By Succession system~\ref{SS} it is clear that any involution
(having $i$ inversions) generates a first involution obtained by
inserting a fixed point of the same parity (it has $i+n-p$
inversions) and a second involution obtained by transforming the
first fixed point by a cycle of a different parity (it has
$i+n-p+1$ inversions).
\\
If $n=2k$ then the number of even involutions in $\II_n(132)$ is
equal to the number of odd involutions in $\II_n(132)$, which is
equal to $|\II_{n-1}(132)|=a_{2k-1}=\binom{2k-1}{k}$.
If $n=4k+1$ then the number of even involutions in $\II_n(132)$ is
equal to $|\II_{n-1}(132)|=a_{4k}=\binom{4k}{2k}$, whereas the
number of odd involutions in $\II_n(132)$ is equal to the number
of involutions in $\II_{n-1}(132)$ excepted those having no fixed
points; that is, $a_{4k}-a_{4k,0}=\binom{4k}{2k+1}$. If $n=4k+3$
then the number of even involutions in $\II_n(132)$ is equal to
the number of involutions in $\II_{n-1}(132)$ excepted those
having no fixed points; that is,
$a_{4k+2}-a_{4k+2,0}=\binom{4k+2}{2k+2}$, whereas the number of
odd involutions in $\II_n(132)$ is equal to
$|\II_{n-1}(132)|=a_{4k+2}=\binom{4k+2}{2k+1}$.
\end{proof}
One can remark that, from Succession system~\ref{SS}, we
immediately have for the number of even involutions in
$\II_n(132)$
$$\sum_{j=0}^{[n/4]} |\{\pi \in \II_n(132) : \pi \mbox{ has }
4j+(n\mod 4) \mbox{ fixed points}\}| = \binom{n-1}{2[(n+1)/4]}$$
and for the number of odd involutions in $\II_n(132)$
$$\sum_{j=0}^{[(n-2)/4]} |\{\pi \in \II_n(132) : \pi \mbox{ has }
4j+((n+2)\mod 4) \mbox{ fixed points}\}| = \binom{n-1}{1+2[(n-2)/4]}.$$

Let $a''(n,p,r,b)$ be the number of involutions $\pi\in\II_n(132)$
having $p$ fixed points and $r$ rises and such that its
decomposition $\pi = \pi' \pi'' x \pi'''$ ($x$ is the first fixed
point and $\pi''$ contains only cycles connected to $\pi'$) leads
to $\pi''=\emptyword$ \resp{$\neq \emptyword$} if and only if
$b=1$ \resp{$0$}. Immediately, we deduce from Succession
system~\ref{SSS} that
$$a''(n,p,r,b)=a''(n-1,p+1-2b,r,1-b)+a''(n-1,p+1-2b,r-1,b)$$
together with $a''(0,0,0,0)=1$.
\begin{theorem}
\label{rises} The number of involutions in $\II_n(132)$ having $r$
rises is given by
$$\binom{[n/2]}{[(r+1)/2]} \binom{[(n-1)/2]}{[r/2]}.$$
\end{theorem}
\begin{proof}
Let $w$ be the bilateral word in bijection with an involution in
$\II_n(132)$ having $r$ rises. So $|w|=n$ and
$|wx|_{xx}+|w|_{\ax\ax}=r$. But $w$ can be regarded as two words
$u$ and $v$ on $\{x,\ax\}^*$ obtained in such a way: $u$
\resp{$v$} consists in all the letters following an $x$
\resp{$\ax$} in $wx$ (taken in the same order). For example, if
$w=\ax xx\ax$ then $u=x\ax$ and $v=xx$. This mapping is
bijective: we can retrieve $w$ from $u$ and $v$ because at each
step we choose either the next letter of $u$ or $v$, and there
is only one way to start. Note also that $|u|=[n/2]$ and
$|v|=[(n+1)/2]$, and that $v$ is ended by $x$. Moreover, if
$n$ is even we have $|u|_x=[(r+1)/2]$ and
$|v|_{\overline{x}}=[r/2]$ otherwise $n$ is odd and we have
$|u|_x=[r/2]$ and $|v|_{\overline{x}}=[(r+1)/2]$ such that
$[n/2]=[(n+1)/2]-1$. So the formula holds in each case.
\end{proof}

\begin{remark}
\label{lrm=rises} An involution in $\II_n(132)$ having $r$ rises
with $0 \leq r < n$ has exactly $n-r$ left-to-right minima (that
is one plus the number of descents).
\end{remark}
\begin{proof}
First, let $e$ be a left-to-right minimum (but not the first one
which is the first element). By definition of a left-to-right
minimum, all the elements left to $e$ are greater than $e$, and so
the element just left to $e$ constitutes a descent.
\\
Secondly, consider a descent in an $132$-avoiding involution $\pi$
that is $e'>e$ and $\pi^{-1}(e')+1=\pi^{-1}(e)$. In order to avoid
$132$, all the elements left to $e'$ must be greater than $e$. So
$e$ is a left-to-right minimum. Moreover, the first element of the
involution is always a left-to-right minimum.
\end{proof}


We complete this subsection by relating these results to
\cite[Theorem~4.6]{GuibertThese} (see also \cite{DGG}). Indeed,
Guibert \cite{GuibertThese} established that all the sets (and
also the others obtained by symmetry operations as inverse, mirror
and complement) $\SS_{n+1}(\mxgddyckI)$, $\SS_{n+1}(\mxgddyckII)$,
$\SS_{n+1}(\mxgddyckIII)$, $\SS_{n+1}(\mxgddyckV)$,
$\SS_{n+1}(\mxgddyckVII)$, $\SS_{n+1}(\mxgddyckVIII)$,
$\SS_{n+1}(\mxgddyckIX)$, $\SS_{n+1}(\mxgddyckX)$,
$\SS_{n+1}(\mxgddyckXI)$, $\SS_{n+1}(\mxgddyckXII)$, and
$\SS_{n+1}(\mxgddyckXIV)$ are in bijection with bilateral words of
length $2n$ enumerated by $a_{2n} = {2n \choose n}$. Thus, by
using the correspondence $\Phi \circ \Xi$, all these sets of
permutations of length $n+1$ avoiding four patterns of length $4$
are in bijection with $\II_{2n}(132)$. Moreover,
Guibert \cite{GuibertThese} also stated that
$|\SS_{n+1}(\mxgddyckXV)|=a_{2n}$.
\subsection{Refinements of bijection $\Psi$}
\label{sec2Psi}

In this subsection, following to \cite{GM}, we study the bijection $\Psi$
between involutions containing $132$ exactly once
(of length $n$ and having almost one fixed point)
and $132$-avoiding involutions
(of length $n-2$ and having the same number of fixed points)
according to one more parameter that is the number of inversions.
So we recall the bijection $\Psi$ given in \cite[Theorem~4.1]{GM}.

Let $\pi\in\II_n$ be an involution containing $132$ exactly once
having $p$ fixed points with $1\leq p\leq n$. We have $\pi = \pi'
x z \pi'' y \pi'''$ with $\pi(x)=x$, $\pi(y)=z$ and $1+x=y<z$ such
that $xzy$ is the only subsequence of type $132$ into $\pi$. In
order to obtain $\sigma$ an involution in $\II_{n-2}(132)$ having
also $p$ fixed points in bijection with $\pi$ by $\Psi$, we
replace the subsequence $xzy$ by a fixed point between $\pi''$ and
$\pi'''$.

\begin{theorem}
\label{Psi-inv} Let $\pi$ be an involution containing $132$
exactly once and let $\sigma$ be the $132$-avoiding involution in
bijection by $\Psi$. If $\pi$ has length $n$, $p$ fixed points
with $1 \leq p \leq n$, and $i$ inversions, then $\sigma$ has
length $n-2$, $p$ fixed points, and $i-2n+2p+3$ inversions.
\end{theorem}
\begin{proof}
Let $\pi = \pi' x z \pi'' y \pi'''$ with $\pi(x)=x$, $\pi(y)=z$
and $1+x=y<z$ and let $\sigma = \sigma' \sigma'' t \sigma'''$ with
$\sigma(t) = t$ and $\sigma(j) \neq j$ for all $1 \leq j < t$ such
that $\sigma',\sigma'',\sigma'''$ corresponds to
$\pi',\pi'',\pi'''$; respectively. Moreover, $\pi'$ contains only
cycles connected to $\pi''$ or $\pi'''$, $\pi''$ contains only
cycles connected to $\pi'$, and $\pi'''$ contains either cycles
connected to $\pi'$ or fixed points; thus, we have
$n-p=2(|\pi'|+1)=2x$.

The number of inversions in $\pi$ but not in $\sigma$ is the sum
of $|\pi'|$ because all the elements of $\pi'$ are greater than
$x$, $|\pi'|-|\pi''|$ because all the elements of $\pi'$ connected
to $\pi'''$ are greater than $z$ (but all the elements of $\pi'$
connected to $\pi''$ are smaller than $z$), $|\pi'|$ because all
the elements of $\pi'$ are greater than $y$, $|\pi''|$ because $x$
is greater than all the elements of $\pi''$ (but $x<z$ and $x<y$),
$|\pi'|-|\pi''|$ because $x$ is greater than all the elements of
$\pi'''$ excepted its fixed points, $|\pi''|+1$ because $z$ is
greater than all the elements of $\pi''$ and than $y$,
$|\pi'|-|\pi''|$ because $z$ is greater than all the elements of
$\pi'''$ excepted its fixed points, $0$ (all the elements of
$\pi''$ are smaller than $y$), and $|\pi'|-|\pi''|$ because $y$ is
greater than all the elements of $\pi'''$ excepted its fixed
points.

The number of inversions in $\sigma$ but not in $\pi$ is the sum
of $|\pi'|-|\pi''|$ because all the elements of $\sigma'$
connected to $\sigma'''$ are greater than $t$ (but all the
elements of $\sigma'$ connected to $\sigma''$ are smaller than
$t$), $0$ (all the elements of $\sigma''$ are smaller than $t$),
and $|\pi'|-|\pi''|$ because $t$ is greater than all the elements
of $\sigma'''$ excepted its fixed points.

Thus, $\sigma$ has
$i-(6|\pi'|+1-2|\pi''|)+(2|\pi'|-2|\pi''|) = i-4|\pi'|-1 = i-2n+2p+3$
inversions.
\end{proof}

\begin{theorem}
\label{Psi-even} The number of even \resp{odd} involutions
containing $132$ exactly once of length $n$ is given by
\begin{center}
$\binom{n-3}{1+2[(n-5)/4]}$\
$\left[{respectively;}\
\binom{n-3}{2[(n-3)/4]}\right]$
\end{center}
for all $n \geq 5$ \resp{$n \geq 3$}.
\end{theorem}
\begin{proof}
By Theorem~\ref{Psi-inv} we trivially deduce that
an even \resp{odd} involution containing $132$ once
corresponds by $\Psi$ to an odd \resp{even} $132$-avoiding involution
because the parity of the number of inversions of the two
involutions is different. Moreover, Theorem~\ref{even132} gives
the number of odd \resp{even} $132$-avoiding involutions according
to the length.

Finally, \cite[Corollary~2.4]{GM} enumerates $132$-avoiding involutions
according to the number of fixed points.
In particular, the number of $132$-avoiding involutions
without fixed points of length $2n$ is $C_n$ the $n$th Catalan number
(indeed $\Phi$ sets these involutions in bijection with Dyck words).
\\
Thus, the number of even \resp{odd} involutions containing $132$
once of length $n$ is $\binom{n-3}{1+2[(n-4)/4]}$
\resp{$\binom{n-3}{2[(n-1)/4]}$} the number of odd \resp{even}
$132$-avoiding involutions of length $n-2$ minus $C_{2j-1}$
\resp{$C_{2j}$} if $n=4j$ \resp{$n=4j+2$} with $j \geq 1$ the
number of $132$-avoiding involutions without fixed points of
(even) length.
\end{proof}
\section{Restricted $132$-even (odd) involutions}
\label{sec3} In this section we study, by block decompositions
approach (see \cite{MV2,MV5}), generating functions for the number
of even or odd involutions on $n$ letters avoiding (or containing
exactly once) $132$ and avoiding (or containing exactly once) an
arbitrary pattern $\tau$ on $k$ letters. The core of this approach
initiated by Mansour and Vainshtein \cite{MV2} lies in the study
of the structure of $132$-avoiding permutations, and permutations
containing a given number of occurrences of $132$.

In several interesting cases the generating function depends only
on $k$ and is expressed via Chebyshev polynomials of the second
kind.

This section is organized in four subsections as corresponding to
the four cases of even or odd involutions avoiding (or containing
exactly once) $132$ and avoiding (or containing exactly once)
$\tau$.
\subsection{Avoiding $132$ and another pattern}

Let $\ie_\tau(n)$ \resp{$\io_\tau(n)$} denote the number of even
\resp{odd} involutions in $\II_n(132,\tau)$, and let
$\ie_\tau(x)=\sum_{n\geq 0}\ie_\tau(n)x^n$
\resp{$\io_\tau(x)=\sum_{n\geq 0}\io_\tau(n)x^n$} be the
corresponding generating function. We denote by $I_\tau(n)$ the
number of involutions in $\II_n(132,\tau)$, and let $I_\tau(x)$ be
the corresponding generating function. In fact, for all $\tau$
\begin{equation}
\label{i} I_\tau(x)=\ie_\tau(x)+\io_\tau(x).
\end{equation}
The following proposition is the base of all the other results
in this section, which holds immediately from definitions.

\begin{proposition}\label{prom1}
Let $\pi\in\II_n(132)$ be an even \resp{odd} involution such that
$\pi_j=n$. Then holds one of the following assertions:
\begin{enumerate}
\item   $\pi_n=n$;

\item   $\pi=(\beta,n,\gamma,\delta,j)$ where
    $1\leq j\leq n/2$, $\delta=\beta^{-1}$ and
    $\beta$ avoids $132$ such that,
    \begin{itemize}
    \item[$(2.1)$] if $j$ even, then $\gamma$ is an $132$-avoiding even \resp{odd} involution
    of length $n-2j$;

    \item[$(2.2)$] if $j$ odd, then $\gamma$ is an $132$-avoiding odd \resp{even} involution
    of length $n-2j$.
    \end{itemize}
\end{enumerate}
\end{proposition}

Our present aim is to find the generating functions $\ie_\tau(x)$
and $\io_\tau(x)$, and since that we need the following lemma.
\begin{lemma}\label{genl}
Let $\{y_n\}_{n\geq0}$ and $\{z_n\}_{n\geq0}$ be two sequences and
the corresponding generating functions we denote by $Y(x)$ and
$Z(x)$; respectively. Then
$$\begin{array}{l}
\sum\limits_{n\geq0}\sum\limits_{j=0}^{[n/4]}y_{2j+1}z_{n-4j}x^n=\dfrac{(Y(x^2)-Y(-x^2))Z(x)}{2x^2};\\
\sum\limits_{n\geq0}\sum\limits_{j=0}^{[n/4]}y_{2j}z_{n-4j}x^n=\dfrac{(Y(x^2)+Y(-x^2))Z(x)}{2}.
\end{array}$$
\end{lemma}

\subsubsection{Pattern $\tau=\varnothing$} The first interesting example is
$\tau=\varnothing$ which presented in \cite{SS}.

\begin{theorem}{\rm(Simion and Schmidt
\cite[Proposition~5]{SS})}\label{tha1} The generating function for
the number of even involutions in $\II_n(132)$ is given by
    $$\ie_\varnothing(x)=\frac{2(1-x)-x^2(C(x^2)-C(-x^2))}{2(1-x+x^2C(-x^2))(1-x-x^2C(x^2))}.$$
\end{theorem}
\begin{proof}
By Proposition~\ref{prom1}, we have exactly three possibilities
for block decomposition of an arbitrary $\pi$ even involution in
$\II_n(132)$. Let us write an equation for $\ie_\varnothing(x)$.
The contribution of the first decomposition above is
$x\ie_\varnothing(x)$. By Lemma~\ref{genl}, the contribution of
the second and the third decomposition above are
$\frac{1}{2}x^2(C(x^2)-C(-x^2))\ie_\varnothing(x)$ and
$\frac{1}{2}x^2(C(x^2)+C(-x^2))\io_\varnothing(x)$; respectively.
Therefore,
$$\ie_\varnothing(x)=1+x\ie_\varnothing(x)
+\frac{1}{2}x^2(C(x^2)-C(-x^2))\ie_\varnothing(x)
+\frac{1}{2}x^2(C(x^2)+C(-x^2))\io_\varnothing(x),$$ where $1$ for
the empty involution. Hence, solving the obtained linear equation
together with \ref{i} we get the desired result.
\end{proof}

\subsubsection{Pattern $\tau=12\dots k$} Let us start by the following
example.
\begin{example}\label{exaa1}
Using Proposition \ref{prom1} we have
    $$I^e_{12}(n)=I^o_{12}(n-2),\quad I^o_{12}(n)=I^e_{12}(n-2).$$
Besides $I^e_{12}(0)=I^e_{12}(1)=1$ and $I^o_{12}(0)=I^o_{12}(1)=0$,
so
$$I^e_{12}(x)=\frac{1+x}{1-x^4},\quad
I^o_{12}(x)=\frac{x^2(1+x)}{1-x^4}.$$
\end{example}

The case of varying $k$ is more interesting. As an extension of
Example \ref{exaa1} let us consider the pattern $\tau=12\dots k$.

\begin{theorem} \label{tha2}
For all $k\geq 1$,
$$\ie_{12\dots k}(x)=
\sum\limits_{j=0}^{k-1}\left(x^j\biggl(1+x^2R_{k-1-j}^e(x^2)I_{k-j}(x)\biggr)\prod\limits_{i=k-j}^k
R_i(-x^2)\right),$$ where $I_{12\dots k}(x)=\frac{1}{x\cdot
U_k\left( \frac{1}{2x}\right)}\cdot\sum\limits_{j=0}^{k-1}
U_j\left( \frac{1}{2x} \right)$ {\rm(see \cite[Theorem~2.8]{GM})}.
\end{theorem}
\begin{proof}
By Proposition \ref{prom1} with the use of the generating function for
the number of permutations in $\SS_n(132, 12\dots k)$ given by
$R_k(x)$ (see \cite{CW}), we have exactly three possibilities for
block decomposition of an arbitrary even involution in
$\II_n(132)$. Let us write an equation for $\ie_{12\dots k} (x)$
with use of Lemma~\ref{genl}. The contribution of the first,
second, and the third decomposition above is
$x\ie_{12\dots(k-1)}(x)$, $x^2\ro_{k-1}(x^2)\ie_{12\dots k}(x)$,
and $x^2\re_{k-1}(x^2)\io_{12\dots k}(x)$; respectively.
Therefore, by definitions of $\re_k(x)$ and $\ro_k(x)$ and by
Identities \ref{eqrk} and \ref{i} we get for all $k\geq 3$,
    $$\ie_{12\dots k}(x)=R_k(-x^2)\left( 1+x^2\re_{k-1}(x^2)I_k(x) \right)
            +xR_k(-x^2)\ie_{12\dots(k-1)}(x).$$
Hence, by the principle of induction on $k$ and by
\cite[Theorem~2.8]{GM} together with Example~\ref{exaa1} we have
the desired result.
\end{proof}

\begin{example}\label{exaa2}
By Theorem~\ref{tha2} and Example~\ref{exaa1} it is easy to see
that
$$I^e_{123}(x)=\frac{1+x+x^2-2x^4}{1-4x^4},$$
and then for all $n\geq 0$,
$$(I^e_{1234}(4n),I^e_{1234}(4n+1),I^e_{1234}(4n+2),I^e_{1234}(4n+3))=(F_{4n-2},F_{4n+1},F_{4n+2},F_{4n+2}),$$
where $F_m$ is the $m$th Fibonacci number.
\end{example}

We can also establish these results for $k=$ $3$, $4$ and $5$ by a
combinatorial approach. In \cite[Subsection~2.1, formula
$(2)$]{GM} it is established that involutions avoiding both $132$
and $1 2 \ldots k$ can be characterized by the succession system:
\begin{equation}
\left\{\begin{array}{lcll}
    (0) \\
    (0) & \leadsto & (1) \\
    (p) & \leadsto & (p+1) , (p-1) & 1 \leq p \leq k-2 \\
    (k-1) & \leadsto & (k-2)
    \end{array}\right.
\label{star}
\end{equation}
and three simple bijections are stated between involutions
avoiding both $132$ and $1 2 \ldots k$ (for $k=3,4,5$) and some
words \resp{$\{a,b\}^*$ or $a\{a,b\}^*$, $\{a,b^2\}^*$,
$\{a,b,c\}^* a$ or $\{a,b,c\}^* a \cup b \{a,b,c\}^* a$}. It is
clear that $|\II_n(132,12\ldots k)|$ for $k=3,4,5$ according to
the number of fixed points $p$ is ($n$ and $p$ have same parity)
is $2^{[(n-1)/2]}$ for $0 \leq p \leq 2$ and $k=3$, $F_{n-2}$
\resp{$F_{n-1}$} for $p=$ $0$ or $3$ \resp{$1$ or $2$} and $k=4$,
and $3^{n/2-1}$ \resp{$\frac{3^{[(n-1)/2]}+1}{2}$, or
$\frac{3^{[(n-1)/2]}-1}{2}$} for $p=2$ \resp{$p=$ $0$ or $1$, or
$p=$ $3$ or $4$} and $k=5$.
\\
So we now deduce from Succession system~\ref{SS} that the number
of even \resp{odd} involutions avoiding both $132$ and $1 2 \ldots
k$ of length $n=4l+1,4l+2,4l+3,4l+4$ with $l \geq 0$ is:
$2^{2l},2^{2l},0,2^{2l+1}$ \resp{$0,2^{2l},2^{2l+1},2^{2l+1}$} for
$k=3$, $F_{4l},F_{4l+1},F_{4l+1},F_{4l+2}$
\resp{$F_{4l-1},F_{4l},F_{4l+2},F_{4l+3}$} for $k=4$, and
$\frac{3^{2l}+1}{2},3^{2l},\frac{3^{2l+1}-1}{2},3^{2l+1}$
\resp{$\frac{3^{2l}-1}{2},3^{2l},\frac{3^{2l+1}+1}{2},3^{2l+1}$}
for $k=5$.
\\
Moreover, we can see this split between even and odd involutions
onto the words in bijection. For $k=3$, the words of $\{a,b\}^n$
or $a\{a,b\}^n$ (in bijection with involutions respectively of
length $2n$ or $2n+1$) corresponds for the even \resp{odd} part to
the words of $a\{a,b\}^{2k} \cup b\{a,b\}^{2k+1}$
\resp{$a\{a,b\}^{2k+1} \cup b\{a,b\}^{2k}$} for any $k \geq 0$.
For $k=4$, the words of $\{a,b^2\}^*$ corresponds for the even
\resp{odd} part to the words of $a\{a,b^2\}^{4k} \cup
a\{a,b^2\}^{4k+1} \cup b^2\{a,b^2\}^{4k+1} \cup
b^2\{a,b^2\}^{4k+2}$ \resp{$a\{a,b^2\}^{4k+2} \cup
a\{a,b^2\}^{4k+3} \cup b^2\{a,b^2\}^{4k+3} \cup
b^2\{a,b^2\}^{4k}$}. We do not give the languages for $k=5$
because their expressions are not so simple.

\subsubsection{Pattern $\tau=2134\dots k$}
Using the argument proof of Theorem~\ref{tha2} together with
$I^e_{21}(x)=\frac{1}{1-x}$ we get

\begin{theorem}\label{tha3}
For all $k\geq 2$,
$$\ie_{2134\dots k}(x)
=\sum\limits_{j=0}^{k-3}\left(x^j\biggl(1+x^2R_{k-1-j}^e(x^2)I_{k-j}(x)\biggr)\prod_{i=k-j}^kR_i(-x^2)\right)
    +x^{k-2}R_2(x)\prod\limits_{j=3}^k R_j(-x^2).$$
where $I_{2134\dots k}(x)=\frac{1}{x\cdot U_k\left(
\frac{1}{2x}\right)}\cdot\sum\limits_{j=0}^{k-1} U_j\left(
\frac{1}{2x} \right)$ {\rm(see \cite[Theorem~2.10]{GM})}.
\end{theorem}

\subsubsection{Pattern $\tau=(d+1,d+2,\dots,k,1,2,\dots,d)$} Let us start
by the following example.

\begin{example}\label{exaa3}
By Proposition~\ref{prom1} it is easy to see that
    $$I^e_{231}(x)=\frac{x^4-x^3+x^2-x+1}{(1-x)^2(1+x^2)}.$$
\end{example}

The case of varying $k$ is more interesting. As an extension of
Example~\ref{exaa3} let us consider the pattern
$\tau=<k,d>:=(d+1,d+2,\dots,k,1,2\dots,d)$.

\begin{theorem}\label{tha4}
For all $1\leq d\leq k/2$,
$$\begin{array}{l}
\ie_{<k,d>}(x)=\frac{R_{k-d}(-x^2)}{1-xR_{k-d}(-x^2)}\Bigl(
1+x^2(R_{k-d-1}(-x^2)-R_{d-1}(-x^2))\ie_{12\dots(k-d)}(x)+\\
\qquad\qquad\qquad\qquad\qquad\qquad\qquad\qquad+(\re_{d-1}(-x^2)-\re_{k-d-1}(-x^2))I_{12\dots(k-d)}(x)+\\
\qquad\qquad\qquad\qquad\qquad\qquad\qquad\qquad\qquad\qquad\qquad\qquad\qquad\qquad+\re_{k-d-1}(x^2)I_{<k,d>}(x)\Bigr),
\end{array}$$
where $I_{12\dots l}(x)$ is given by \cite[Theorem~2.10]{GM} and
$I_{<k,d>}(x)$ is given by \cite[Theorem~2.14]{GM}.
\end{theorem}
\begin{proof}
By Proposition \ref{prom1}, we have exactly three possibilities
for block decomposition an arbitrary even involution
$\pi\in\II_n(132)$. Let us write an equation for $\ie_{<k,d>}(x)$
by using Lemma~\ref{genl}. The contribution of the first above
decomposition above is $x I^e_{<k,d>}(x)$. To find the
contribution of the second and the third decompositions above let
us consider two cases. First, if $\gamma$ avoids $12\dots (k-d)$,
then $\beta$ and $\delta$ avoid $12\dots (k-d-1)$, so the
generating function for this number of even involutions is given
by
$$x^2\ro_{k-d-1}(x^2)\ie_{12\dots(k-d)}(x)+x^2\re_{k-d-1}(x^2)\io_{12\dots(k-d)}.$$
Secondly, if $\gamma$ contains $12\dots (k-d)$,
then $\beta$ and $\delta$ avoid $12\dots (d-1)$, so the
generating function for this number of even involutions is given
by
$$x^2\ro_{d-1}(x^2)(\ie_{<k,d>}(x)-\ie_{12\dots(k-d)}(x))+
x^2\re_{d-1}(x^2)(\io_{<k,d>}(x)-\io_{12\dots(k-d)}(x)).$$

We mentioned that, the generating function for the number of even
involutions in $\II_n(132,<k,d>)$ such containing $12\dots(k-d)$
at least once is given $\ie_{<k,d>}(x)-\ie_{12\dots(k-d)}(x)$.
Therefore,
$$\begin{array}{l}
\ie_{<k,d>}(x)=1+x\ie_{<k,d>}(x)+x^2\ro_{k-d-1}(x^2)\ie_{12\dots(k-d)}(x)+x^2\re_{k-d-1}(x^2)\io_{12\dots(k-d)}(x)+\\
\qquad\qquad\qquad+x^2\ro_{d-1}(x^2)(\ie_{<k,d>}(x)-\ie_{12\dots(k-d)}(x))+x^2\re_{d-1}(x^2)(\io_{<k,d>}(x)-\io_{12\dots(k-d)}(x)),
\end{array}$$
Hence, by Identities \ref{eqrk}, \ref{eqeo} and \ref{i} we get the
desired result.
\end{proof}

%
\subsection{Avoiding $132$ and containing another pattern}

Let $\ie_{\tau;r}(n)$ \resp{$\io_{\tau;r}(n)$} denote the number
of even \resp{odd} involutions in $\II_n(132)$ such containing
$\tau$ exactly $r$ times, and let
$\ie_{\tau;r}(x)=\sum_{n\geq0}\ie_{\tau;r}(n)x^n$ \resp{
$\io_{\tau;r}(x)=\sum_{n\geq0}\io_{\tau;r}(n)x^n$} be the
corresponding generating function. Also, we denote the generating
function for the number of $132$-avoiding involutions containing
$\tau$ exactly $r$ times by $I_{\tau;r}(x)$. So for all $\tau$
\begin{equation}
\label{ii}
I_{\tau;r}(x)=\ie_{\tau;r}(x)+\io_{\tau;r}(x).
\end{equation}

\subsubsection{Pattern $\tau=12\dots k$}
Now let us start by the following example.

\begin{example}\label{exbb1}
Proposition \ref{prom1} yields
    $$\ie_{12;1}(x)=x^2+x^2\io_{12;1}(x),\quad \io_{12;1}(x)=x^2\ie_{12;1}(x),$$
which means that $\ie_{12;1}(x)=\frac{x^2}{1-x^4}$ and
$\io_{12;1}(x)=\frac{x^4}{1-x^4}$.
\end{example}

The case of varying $k$ is more interesting. As an extension of
Example~\ref{exbb1} let us consider the pattern $\tau=12\dots k$.

\begin{theorem} \label{thb2}
For all $k\geq 2$,
    $$\ie_{12\dots k;1}(x)=\sum_{j=0}^{k-1}
    \frac{x^{k-j}R_{j+1}^e(x^2)\prod_{i=j+2}^k R_i(-x^2)}
    {U_{j+2}\left(\frac{1}{2x}\right)}.$$
\end{theorem}
\begin{proof}
By Proposition~\ref{prom1}, we have exactly three possibilities
for block decompositions an arbitrary even involution in
$\II_n(132)$. Similarly as proof of Theorem~\ref{tha3} we have
$$\ie_{12\dots
k;1}(x)=x\ie_{12\dots(k-1);1}(x)+x^2\ro_{k-1}(x^2)\ie_{12\dots
k;1}(x)+x^2\re_{k-1}(x^2)\io_{12\dots k;1}(x).$$ hence, by using
Identities \ref{eqrk}, \ref{eqeo}, and \ref{ii} together with
\cite[Theorem~3.2]{GM}, and then using the principle of induction
on $k$ together with Example~\ref{exbb1} we get the desired
result.
\end{proof}

\begin{example}\label{exbb2}
Theorem \ref{thb2} yields
$\ie_{1234;1}(x)=\frac{x^4(x^4+1)}{(1+3x^2+x^4)(1-3x^2+x^4)}$. In
other words, for all $n\geq 0$ we have $\ie_{1234;1}(n)=F_{n-2}$
if $n/4$ is a positive integer number otherwise
$\ie_{1234;1}(n)=0$, where $F_m$ is the $m$th Fibonacci number.
\end{example}

\subsubsection{Pattern $\tau=2134\dots k$} Similarly as Theorem~\ref{tha2}
and Theorem~\ref{thb2} together with using \cite[Theorem~3.3]{GM}
we have the following result.
\begin{theorem}\label{thb3}
For all $k\geq 2$,
    $$\ie_{2134\dots k;1}(x)=(1-x^2)\sum_{j=3}^k
    \frac{x^{k+2-j}\re_{j-1}(x^2)\prod_{i=j}^kR_i(-x^2)}
    {U_j\left(\frac{1}{2x}\right)}.$$
\end{theorem}

\begin{example}\label{exbb2bis}
Theorem~\ref{thb2} yields
$\ie_{2134;1}(x)=\frac{x^6(2-x^4)}{(1+3x^2+x^4)(1-3x^2+x^4)}$. In
other words, for all $n\geq 0$ we have $\ie_{2134;1}(n)=F_{n-3}$
if $n/4$ is a positive integer number otherwise
$\ie_{2134;1}(n)=0$, where $F_m$ is the $m$th Fibonacci number.
\end{example}

\subsubsection{Pattern $\tau=23\dots k1$}
Now, let us consider another interesting case where $\tau=23\dots
k1$.

\begin{theorem}\label{thb4}
For all $k\geq 3$,
$$\ie_{23\dots k1;1}(x)=\frac{x^3}{1-x}\io_{12\dots(k-2);1}(x).$$
\end{theorem}
\begin{proof}
By Proposition~\ref{prom1} it is easy to see that
$$\ie_{23\dots k1;1}(x)=x\ie_{23\dots k1;1}(x)+x^2g(x),$$
where $g(x)$ the generating function for the number of odd
involutions in $\II_n(132,23\dots k1)$ such containing
$12\dots(k-1)$ exactly once. On the other hand, also by
Proposition~\ref{prom1} we get
    $$g(x)=x\io_{12\dots(k-2);1}(x).$$
Hence the theorem holds.
\end{proof}

As a remark, to find an explicit formula for $\ie_{23\dots
k1;1}(x)$ see Theorem~\ref{thb4}, Theorem~\ref{thb2}, and
Identity~\ref{ii}.

\begin{example}\label{exbb3}
Theorem~\ref{thb4} and Theorem~\ref{thb2} yield $\ie_{231;1}(x)=0$
and $\ie_{2341;1}(x)=\frac{x^7}{(1-x)(1-x^4)}$.
\end{example}

%
\subsection{Containing $132$ once and avoiding another pattern}

Let $\je_\tau(n)$ \resp{$\jo_\tau(n)$} denote the number of even
\resp{odd} involutions in $\II_n(\tau)$ such containing $132$
exactly once, let $\je_\tau(x)=\sum_{n\geq0}\je_\tau(n)x^n$ \resp{
$\jo_\tau(x)=\sum_{n\geq0}\jo_\tau(n)x^n$} be the corresponding
generating function. Also, we denote by $J_\tau(x)$ the generating
function for the number of involutions in $\II_n(\tau)$ containing
$132$ exactly $r$ times. So we have
\begin{equation}
J_\tau(x)=\je_\tau(x)+\jo_\tau(x).
\label{iii}
\end{equation}

The following proposition is the base of all the other results
in this section, which holds immediately from definitions.

\begin{proposition}\label{prom2}
Let $\pi\in\II_n$ be an even involution contains $132$ exactly
once such that $\pi_j=n$. Then holds one of the following
assertions:
\begin{enumerate}
\item   $\pi_n=n$;
\item   $\pi=(\beta,n,\gamma,\delta,j)$
    where $1\leq j\leq n/2$, $\delta=\beta^{-1}$, and
    $\beta$ avoids $132$, such that
    \begin{itemize}
    \item[$(2.1)$] if $j$ even, then $\gamma$ is an even involution contains $132$ exactly once
    of length $n-2j$,
    \item[$(2.2)$] if $j$ odd, then $\gamma$ is an odd involution contains $132$ exactly once
    of length $n-2j$;
    \end{itemize}
\item   $\pi=(\beta,m,2m+1,\gamma,m+1)$ where $n=2m+1$ and
    $\gamma=\beta^{-1}$ such that
    $m$ even and $\beta\in\SS_{m-1}(132)$.
\end{enumerate}
\end{proposition}

\subsubsection{Pattern $\tau=\varnothing$} The first interesting case is where
$\tau=\varnothing$ which is analogue to Theorem~\ref{tha1}.

\begin{theorem}\label{thc1} The generating function for the number
of involutions in $\II_n$ containing $132$ exactly once is given
by
$$\je_\varnothing(x)=\frac{x\left( 1-2x+x\sqrt{1-4x^2}+x\sqrt{1+4x^2}-\sqrt{1-4x^2}\sqrt{1+4x^2}\right)}
    {(1-2x+\sqrt{1+4x^2})(1-2x+\sqrt{1-4x^2})}.$$
In other words, for all $n\geq 1$,
    $$\je_\varnothing(n)=\frac{1}{2}C_{(n-2)/2}\left( \frac{1}{2}[n/2](3+(-1)^{n+1})-1-(-1)^{\binom{n-3}{2}} \right),$$
where $C_m$ is the $m$th Catalan number.
\end{theorem}
\begin{proof}
Proposition \ref{prom2}, we have exactly four possibilities for
block decomposition an arbitrary even involution in $\II_n$
containing $132$ exactly once. Let us write an equation for
$J^e_\varnothing(x)$ by using Lemma~\ref{genl}. The contribution
of the first, the second, the third, and the fourth decompositions
above are $x\je_\varnothing(x)$,
$\frac{x^2}{2}(C(x^2)-C(-x^2))\je_\varnothing(x)$,
$\frac{x^2}{2}(C(x^2)+C(-x^2))\jo_\varnothing(x)$, and
$\frac{x^3}{2}(C(x^2)-C(-x^2))$; respectively. Therefore,
$$\je_\varnothing(x) = x\je_\varnothing(x) +
\frac{x^2}{2}(C(x^2)-C(-x^2))\je_\varnothing(x) +
\frac{x^2}{2}(C(x^2)+C(-x^2))\jo_\varnothing(x) +
\frac{x^3}{2}(C(x^2)-C(-x^2)).$$
Hence, by Identity~\ref{iii} and \cite[Theorem~4.4]{GM} we get the
desired result.
\end{proof}

\subsubsection{Pattern $\tau=12\dots k$} Let us start by the following
example.

\begin{example}\label{excc1}
By Proposition~\ref{prom2} it is easy to see $\je_{12}(x)=0$.
\end{example}

As extension of Example~\ref{excc1} let us consider the pattern
$\tau=12\dots k$.

\begin{theorem}\label{thc2}
For all $k\geq 2$,
$$\je_{12\dots k}(x)=\sum\limits_{j=3}^k\left(
x^{k+1-j} \biggl[\ro_{j-1}(x^2)+
\frac{\re_{j-1}(x^2)}{U_j\left(\frac{1}{2x}\right)}
  \sum\limits_{i=1}^{j-2}U_i(1/(2x))\biggr]\prod\limits_{i=j}^k R_i(-x^2)\right).$$
\end{theorem}
\begin{proof}
Similarly as the argument proof of Theorem~\ref{thc1} together
with using the fact that the generating function for the number of
permutations in $\SS_n(132,12\dots m)$ is given by $R_m(x)$ (see
\cite{CW}), we get
$$\je_{12\dots k}(x)=x\je_{12\dots
(k-1)}(x)+x^2\ro_{k-1}(x^2)\je_{12\dots
k}(x)+x^2\re_{k-1}(x^2)\jo_{12\dots k}(x)+x^3\ro_{k-1}(x^2).$$
Therefore, by Identities \ref{eqrk} and \ref{iii} we have
$$\je_{12\dots k}(x)=xR_k(-x^2)(x^2\ro_{k-1}(x^2)+x\re_{k-1}(x^2)J_{12\dots
k}(x)+\je_{12\dots(k-1)}(x)).$$ Hence, by the principle of
induction on $k$ together with Example~\ref{excc1} and
\cite[Theorem~4.6]{GM} we get the desired result.
\end{proof}

\begin{example}\label{excc2}
Theorem~\ref{thc2} for $k=4$ yields $\je_{1234}(n)=F_{n-2}$ where
$n=5,6,9,10,13,14,\dots$, otherwise $\je_{1234}(n)=0$, where
$F_{m}$ is the $m$th Fibonacci number.
\end{example}

\subsubsection{Pattern $\tau=2134\ldots k$ or $\tau=23\ldots k1$}
Similarly as the argument proof of Theorem~\ref{thc2} together
with using the fact that the generating function for
$132$-avoiding permutations such avoiding $2134\ldots k$ (or
$23\ldots k1$) is given by $R_k(x)$ (see \cite{MV2}), and use
\cite[Theorems~4.7 and 4.10]{GM} (expressions for $J_{2134\ldots
k}(x)$ and $J_{23\ldots k1}(x)$; respectively) we obtain other
cases $\tau=2134\ldots k$ and $\tau=23\ldots k1$.

\begin{theorem}\label{thc3}
$$\je_{2134\dots k}(x) =xR_k(-x^2)\Bigl(x^2\ro_{k-1}(x^2)+x\re_{k-1}(x^2)
J_{2134\dots k}(x)+\je_{2134\dots(k-1)}(x)\Bigr),$$ with
$\je_{213}(x)=\frac{x^5}{1-4x^4}$, where $J_{2134\dots
m}(x)=\frac{x\left[ x^2U_2\left(
\frac{1}{2x}\right)+\sum_{j=2}^{m-2}U_j\left(\frac{1}{2x}\right)\right]}
{U_m\left( \frac{1}{2x}\right) }$ {\rm(}see
\cite[Theorem~4.7]{GM}{\rm)}.
\end{theorem}


\begin{theorem}\label{thc4}
For all $k\geq 2$,
$$\je_{23\dots k1}(x)=\frac{x^2}{1-x}\Bigl[ x\ro_{k-2}(x^2)+\re_{k-2}(x^2)J_{12\dots(k-1)}(x)
-R_{k-2}(-x^2)\je_{12\dots(k-1)}(x) \Bigr],$$ where $J_{12\dots
m}(x)=\frac{x}{U_m\left( \frac{1}{2x} \right)} \sum_{j=1}^{m-2}
U_{j}\left( \frac{1}{2x} \right)$ {\rm(}see
\cite[Theorem~4.6]{GM}{\rm)}.
\end{theorem}

%
\subsection{Containing $132$ and another pattern once}

Let $\je_{\tau;r}(n)$ \resp{$\jo_{\tau;r}(n)$} denote the number
of even \resp{odd} involutions in $\II_n$ such containing $132$
exactly once and containing $\tau$ exactly $r$ times, let
$\je_{\tau;r}(x)=\sum_{n\geq0}\je_{\tau;r}(n)x^n$ \resp{
$\jo_{\tau;r}(x)=\sum_{n\geq0}\jo_{\tau;r}(n)x^n$} be the
corresponding generating function.

\begin{example}\label{exdd1}
By Proposition~\ref{prom2} it is easy to see that
$\je_{12;1}(x)=0$ and $\je_{21;1}(x)=0$. So the first interesting
case can be examined in this case is when $\tau=123$, $231$, or
$213$. But $\je_{123;1}(x)=\je_{231;1}(x)=0$ and
$\je_{213;1}(x)=\frac{x^6}{1-4x^4}$.
\end{example}

Once again, as extension of Example~\ref{exdd1} let us consider
the patterns $\tau=12\dots k$, $\tau=23\dots k1$, and $\tau=2134\dots
k$. Similarly as the arguments proofs in the last subsection
together with \cite[Theorems~5.1, 5.2 and 5.4]{GM} we obtain as
the following.

\begin{theorem}\label{thd1}
For all $k\geq 2$,
        $$\je_{12\dots k;1}(x)=\je_{23\dots k1;1}(x)=0.$$
\end{theorem}

\begin{theorem}\label{thd2}
For all $k\geq 2$,
    $$\je_{2134\dots k;1}(x)=(1-x^2)\sum_{j=2}^{k-1}
    \frac{x^{k+2-j}\re_j(x^2)\prod_{i=j+1}^kR_i(-x^2)}
    {U_{j+1}\left(\frac{1}{2x}\right)}.$$
\end{theorem}

%

\section{Statistics of generalized patterns on $132$-avoiding involutions}
\label{sec4}

In the current section let us consider the case of generalized
patterns (see Subsection~\ref{intro-back} for their definitions),
and let us study some statistics on $1\mn3\mn2$-avoiding involutions.
We relate
some of these results to the enumeration of $1\mn3\mn2$-avoiding
involutions according to the length and the number of rises (given
in Subsection~\ref{sec2Phi}).

Robertson, Wilf and Zeilberger \cite{RWZ} showed a simple
continued fraction that records the joint distribution of the
patterns $1\mn2$ and $1\mn2\mn3$ on $1\mn3\mn2$-avoiding
permutations. Recently, generalization of this theorem given, by
Mansour and Vainshtein \cite{MV1}, by Krattenthaler \cite{Kr}, by
Jani and Rieper \cite{JR}, and by Br\"and\"en, Claesson and
Steingr\`{\i}msson \cite{BCS}. Mansour \cite{MS} generalizes the main
result of \cite{BCS} by replacing generalized patterns with
classical patterns.

An another analogue of these results is to replace
$1\mn3\mn2$-avoiding permutations with $1\mn3\mn2$-avoiding
involutions.

First of all, let us define the set of all $1\mn3\mn2$-avoiding
involutions of all sizes with the empty permutation by
$\mathcal{I}(1\mn3\mn2)$, and the set of all $1\mn3\mn2$-avoiding
permutations of all sizes with the empty permutation by
$\mathcal{S}(1\mn3\mn2)$. Now, let us start by following
proposition \cite[Proposition~2.1]{GM} which is the base of all
the other results in the following subsections.

\begin{proposition}\label{prom}{\rm(Guibert and Mansour
\cite[Proposition~2.1]{GM})} Let $\pi\in\II_n(1\mn3\mn2)$ be an
involution such that $\pi_j=n$. Then holds one of the following
assertions:
\begin{enumerate}
\item   for $1\leq j\leq [n/2]$, $\pi=(\beta,n,\gamma,\delta,j)$, where
    \begin{itemize}
    \item[1.1] $\beta$ is a $1\mn3\mn2$-avoiding permutation of the numbers $n-j+1,\dots,n-2,n-1$,

    \item[1.2] $\delta$ is a $1\mn3\mn2$-avoiding permutation of the numbers $1,\dots,j-2,j-1$ such
    that $\delta\cdot\beta$ is the identity permutation of $\SS_{j-1}$,

    \item[1.3] $\gamma$ is a $1\mn3\mn2$-avoiding involution of the numbers $j+1$, $j+2$, $\dots$, $n-j-1$,
    $n-j$;
    \end{itemize}
\item   for $j=n$, $\pi=(\beta,n)$ where $\beta$ is an involution in $\II_{n-1}(1\mn3\mn2)$.
\end{enumerate}
\end{proposition}
\subsection{Counting an occurrences of $1\mn2\mn3\dots\mn k$}

Let us define $$\begin{array}{l}
C_I(x_1,x_2,\dots)=\sum\limits_{\pi\in\mathcal{I}(1\mn3\mn2)}
\prod\limits_{k\geq 1} x_k^{\sharp1\mn2\mn\dots\mn k(\pi)},\\
C_S(x_1,x_2,\dots)=\sum\limits_{\pi\in\mathcal{S}(1\mn3\mn2)}\prod\limits_{k\geq
1} x_k^{\sharp1\mn2\mn\dots\mn k(\pi)}, \end{array}$$ where
$\sharp\tau(\pi)$ is the number of occurrences of $\tau$ in $\pi$.

\begin{theorem}\label{thc}
The generating function $C_I(x_1,x_2,\dots)$ given by
$$C_I(x_1,x_2,\dots)=\dfrac{1+x_1C_I(x_1x_2,x_2x_3,\dots)}
                          {1-x_1^2C_S(x_1^2x_2^2,x_2^2x_3^2,\dots)};$$
where {\rm(}see \cite[Theorem~1]{BCS}{\rm)}
$$C_S(x_1,x_2,\dots)=\dfrac{1}{1-x_1C_S(x_1x_2,x_2x_3,\dots)}.$$
\end{theorem}
\begin{proof}
In \cite[Theorem~1]{BCS} proved
$$C_S(x_1,x_2,\dots)=\frac{1}{1-x_1C_S(x_1x_2,x_2x_3,\dots)}.$$
On the other hand, by Proposition~\ref{prom} it is easy to see for
$k\geq 1$,
$$\sharp1\mn2\mn3\mn\dots\mn k(\pi',n)=\sharp1\mn2\mn3\mn\dots\mn
k(\pi')+\sharp1\mn2\mn3\mn\dots\mn(k-1)(\pi'),$$  and
$$\sharp1\mn2\mn3\mn\dots\mn k(\pi',n,\beta,\pi'',j)=
2\cdot\sharp1\mn2\mn3\mn\dots\mn k(\pi')
+2\cdot\sharp1\mn2\mn3\mn\dots\mn (k-1)(\pi')+
\sharp1\mn2\mn3\mn\dots\mn k(\beta).$$ Let us write an equation
for the generating function $C_I(x_1,x_2,\dots)$. The contribution
of the first case is $x_1C_I(x_1x_2,x_2x_3,\dots)$, and of the
second case is
$x_1^2C_I(x_1,x_2,\dots)C_S(x_1^2x_2^2,x_2^2x_3^2,\dots)$ (see
Lemma~\ref{genl}). Hence
$$C_I(x_1,x_2,\dots)=1+x_1C_I(x_1x_2,x_2x_3,\dots)+x_1^2C_I(x_1,x_2,\dots)C_S(x_1^2x_2^2,x_2^2x_3^2,\dots),$$
where $1$ for empty involution.
\end{proof}

\begin{example}
It is easy to see by Theorem~\ref{thc} that, the number of
involutions in $\II_n(1\mn3\mn2)$ containing $1\mn2$ exactly once
is given by $\frac{1}{2}(1+(-1)^n)$, and the number of involutions
in $\II_n(1\mn3\mn2)$ containing $1\mn2$ exactly twice is given by
$\frac{1}{4}(2n-3-(-1)^n)$.
\end{example}
An application for Theorem~\ref{thc} we get
$$\sum\limits_{\pi\in\mathcal{I}(1\mn3\mn2)}x^{|\pi|}=
\frac{1}{1-x-\dfrac{x^2}{1-\dfrac{x^2}{1-\dfrac{x^2}{\ddots}}}} =
\frac{2}{1-2x+\sqrt{1-4x^2}},$$ which means that the number of
$1\mn3\mn2$-avoiding involutions in $\II_n$ is given by the
central binomial coefficient $a_n=\binom{n}{[n/2]}$ (see \cite{SS}).

Let $s_\pi$ be the number of right-to-left maxima of
$\pi\in\mathcal{S}(132)$; so by \cite[Proposition~5]{BCS} we get
$$s_\pi=\sharp1(\pi)-\sharp1\mn2(\pi)+\sharp1\mn2\mn3(\pi)-\cdots.$$
An application for Theorem~\ref{thc} with \cite[Theorem~1]{BCS}
($C_S(x,1,1,\dots)=C(x)$) we get
$$\sum_{\pi\in\mathcal{I}(1\mn3\mn2)}x^{|\pi|}y^{s_\pi}=\frac{1+xyC_I(x,1,1,\dots)}{1-x^2y^2C_S(x^2,1,1,\dots)},$$
and $C_I(x,1,1,\dots)=\frac{1}{1-x-x^2C(x^2)}$. Hence,
$$\sum_{\pi\in\mathcal{I}(1\mn3\mn2)}x^{|\pi|}y^{s_\pi}
=\sum_{j\geq0}x^{2j}C^j(x^2)y^{2j}+\sum_{j\geq0}\frac{x^{2j+1}C^j(x^2)}{1-x-x^2C(x^2)}y^{2j+1}.$$

We can also prove this last result
$\sum_{\pi\in\mathcal{I}(1\mn3\mn2)}x^{|\pi|}y^{s_\pi}$ from an
entire combinatorial way (for convenience, we put $|\pi|=n$ and
$s_\pi=s$).

\begin{theorem}
\label{th-ssd} The number of $1\mn3\mn2$-avoiding involutions of
length $n$ having $s$ right-to-left maxima with $1 \leq s \leq n$
is:
$$\begin{array}{ll}
0
    & \mbox{if } n \mbox{ is odd and } s \mbox{ is even,} \\
{n-1-s/2 \choose n/2-1} - {n-1-s/2 \choose n/2}
    & \mbox{if } n \mbox{ and } s \mbox{ are even,} \\
{n-1-(s-1)/2 \choose {[n/2]}}
    & \mbox{if } s \mbox{ is odd.} \\
\end{array}$$
\end{theorem}

\begin{remark}
\label{rmq-ssd} Let $\pi$ be an $1\mn3\mn2$-avoiding involution of
length $n$ having $s$ right-to-left maxima (with $1 \leq s \leq
n$) namely $m_1,m_2,\ldots,m_s$ with
$0=m_0<m_1<m_2<\cdots<m_{s-1}<m_s=n$. Thus, we have the following
facts.
\begin{itemize}
\item[(1)]
All the elements of $]m_{i-1},m_i[$ are located between $m_{i+1}$ and $m_i$.
\item[(2)]
$\pi^{-1}(m_i) = m_{s+1-i}$ for all $1 \leq i \leq s$.
\item[(3)]
If $s$ is even, then $\pi$ has no fixed point, else, and if $s
\neq 1$, we have that $\pi$ has all its fixed points between
$m_{\frac{s+1}{2}+1}$ (excluded) and $m_{\frac{s+1}{2}}$ (included
because it is a fixed point), and that
$m_{\frac{s+1}{2}}+m_{\frac{s+1}{2}-1}=n$.
\end{itemize}
\end{remark}
\begin{proof}
(1) is immediately deduced from the definition of a right-to-left
maximum and from the pattern $1\mn3\mn2$ we must avoid. Thus
$\pi^{-1}(m_i) = | ]m_{i-1},n] | = n - m_{i-1}$ for all $1 \leq i
\leq s$.
\\
(2) is deduced from (1), from that $\pi^{-1}(n) = m_1$ and because
$\pi$ is an involution. The proof is successively stated for
$i=1,2,\ldots,[(s+1)/2]$.
\\
(3) is deduced from (1), (2) and Proposition~\ref{prom}.
The last fact is given by a simple calculus:
$m_{\frac{s+1}{2}-1} = \pi^{-1}(m_{\frac{s+1}{2}+1}) = n - m_{\frac{s+1}{2}}$.
\end{proof}

\begin{lemma}
\label{lemm-ssd1}
There is a bijection between $1\mn3\mn2$-avoiding involutions
of length $2k$ having $2l+1$ right-to-left maxima and
$1\mn3\mn2$-avoiding involutions of length $2k+1$ having $2l+3$
right-to-left maxima, with $0 \leq l < k$.
\end{lemma}

\begin{proof}
Let $\pi \in \II_{2k}(1\mn3\mn2)$ having $2l+1$ right-to-left maxima
which are $m_1,m_2,\ldots,m_{2l+1}$ with
$0=m_0<m_1<m_2<\cdots<m_{2l}<m_{2l+1}=2k$. $\pi$ has two or more
fixed points because one of them is a right-to-left maximum by (3)
of Remark~\ref{rmq-ssd}, and the length is even. Let $x$ be the
penultimate fixed point of $\pi$ (the last one is $m_{l+1}$). In
order to avoid $1\mn3\mn2$ and by (1) of Remark~\ref{rmq-ssd}, we
have that the elements of $\pi$ located between $m_{l+2}$ and
$m_{l+1}$ forms a factor $\pi' \pi'' x \pi'''$ with
$|\pi'|=|\pi'''|=m_{l+1}-1-x$ and all the elements of $\pi'$,
$\pi''$ and $\pi'''$ belong respectively to $]x,m_{l+1}[$,
$[m_{l+1}+m_l-x,x[$ and $]m_l,m_{l+1}+m_l-x[$ that is all the
elements of $\pi'$ are connected to all the elements of $\pi'''$.
Thus, we obtain $\sigma$ by changing into $\pi$ the fixed point
$m_{l+1}$ by a cycle starting between $\pi'$ and $\pi''$. $\sigma$
avoids $1\mn3\mn2$, has length $2k+1$ and has $2l+3$ right-to-left
maxima which are
    $m_1$, $m_2$, $\ldots$, $m_l$,
    $m_{l+1}+m_l-x$ (instead of $m_{l+1}$),
    $x+1$ (new one),
    $m_{l+1}+1$ (new one),
    $m_{l+2}+1$, $m_{l+3}+1$, $\ldots$, $m_{2l+1}+1$.
This mapping is clearly bijective.
\end{proof}

For example,
$\pi = 19 \ 18 \ 20 \ 17 \ 15 \ 13 \ 14 \ 10 \ 9 \ 8 \ 11 \ 12 \ 6 \ 7 \ 5 \ 16 \ 4 \ 2 \ 1 \ 3$
of length $20$ ($k=10$) having $5$ right-to-left maxima
($l=2$ and $m_1=3,m_2=4,m_3=16,m_4=17,m_5=20$) with $x=12$,
$\pi' = 15 \ 13 \ 14$, $\pi'' = 10 \ 9 \ 8 \ 11$, $\pi''' = 6 \ 7 \ 5$,
and $\sigma = 20 \ 19 \ 21 \ 18 \ 16 \ 14 \ 15 \ 17 \ 11 \ 10 \ 9 \ 12 \ 13 \ 6 \ 7 \ 5 \ 8 \ 4 \ 2 \ 1 \ 3$
of length $21$ having $7$ right-to-left maxima ($3,4,8,13,17,18,21$)
are in bijection.

\begin{lemma}
\label{lemm-ssd2}
There is a bijection between
$1\mn3\mn2$-avoiding involutions
    of length $2k+1$ having $2l+1$ right-to-left maxima
and $1\mn3\mn2$-avoiding involutions
    of length $2k+2$ having $2l+2$ or $2l+3$ right-to-left maxima,
with $0 \leq l \leq k$.
\end{lemma}

\begin{proof}
Let $\pi \in \II_{2k+1}(1\mn3\mn2)$
having $2l+1$ right-to-left maxima  which are $m_1,m_2,\ldots,m_{2l+1}$
with $0=m_0<m_1<m_2<\cdots<m_{2l}<m_{2l+1}=2k+1$.
There are two cases to consider depending on the number
of fixed points of $\pi$.

The first case consists in considering $\pi$ having one fixed point
which is $m_{l+1}$. In order to avoid $1\mn3\mn2$ and because
there is only one fixed point, we have that the elements of $\pi$
located between $m_{l+2}$ and $m_{l+1}$ forms a factor $\pi'
\pi''$ with $|\pi'|=|\pi''|$ such that all the elements of $\pi'$
are connected to all the elements of $\pi''$. Thus, $\sigma$ the
involution of even length is obtained by changing into $\pi$ the
fixed point $m_{l+1}$ by a cycle starting between $\pi'$ and
$\pi''$. $\sigma$ avoids $1\mn3\mn2$, has length $2k+2$, has no
fixed point and has $2l+2$ right-to-left maxima which are
    $m_1$, $m_2$, $\ldots$, $m_l$,
    $m_{l+1}+m_l$ (instead of $m_{l+1}$),
    $m_{l+1}+1$ (new one),
    $m_{l+2}+1$, $m_{l+3}+1$, $\ldots$, $m_{2l+1}+1$.
This mapping is clearly bijective.

The second case consists in considering $\pi$ having
more than two fixed points, and we apply exactly the same
bijection given in Lemma~\ref{lemm-ssd1}.
The involution $\sigma$ we obtain avoids $1\mn3\mn2$, has length
$2k+2$ and has $2l+3$ right-to-left maxima.
\end{proof}

As an example of the first case, $\pi = 18 \ 17 \
19 \ 14 \ 15 \ 16 \ 11 \ 10 \ 12 \ 8 \ 7 \ 9 \ 13 \ 4 \ 5 \ 6 \ 2
\ 1 \ 3$ of length $19$ ($k=9$) having $5$ right-to-left maxima
($l=2$ and $m_1=3,m_2=6,m_3=13,m_4=16,m_5=19$) with $\pi' = 11 \
10 \ 12$, $\pi'' = 8 \ 7 \ 9$ and one fixed point ($13$), and
$\sigma = 19 \ 18 \ 20 \ 15 \ 16 \ 17 \ 12 \ 11 \ 13 \ 14 \ 8 \ 7
\ 9 \ 10 \ 4 \ 5 \ 6 \ 2 \ 1 \ 3$ of length $20$ having $6$
right-to-left maxima ($3,6,10,14,17,20$) and no fixed point are in
bijection.
\\
Moreover, as an example of the second case, $\pi = 12
\ 13 \ 10 \ 9 \ 6 \ 5 \ 7 \ 8 \ 4 \ 3 \ 11 \ 1 \ 2$ of length $13$
($k=6$) having $3$ right-to-left maxima ($l=1$ and
$m_1=2,m_2=11,m_3=13$) with $x=8$, $\pi' = 10 \ 9$, $\pi'' = 6 \ 5
\ 7$, $\pi''' = 4 \ 3$ and $3$ fixed points ($7,8,11$) and $\sigma
= 13 \ 14 \ 11 \ 10 \ 12 \ 7 \ 6 \ 8 \ 9 \ 4 \ 3 \ 5 \ 1 \ 2$ of
length $14$ having $5$ right-to-left maxima ($2,5,9,12,14$) and
$2$ fixed points ($8,9$) are in bijection.

{\em Proof of Theorem~\ref{th-ssd}}. Let $\pi\in\II_n(1\mn3\mn2)$
having $s$ right-to-left maxima with $1 \leq s \leq n$.
\\
By $(3)$ of Remark~\ref{rmq-ssd} we immediately obtain that there is no
$\pi$ such that $n$ is odd and $s$ is even. Trivially, if $\pi$
has only one right-to-left maximum ($s=1$), the $n-1$ first
elements constitute an $1\mn3\mn2$-avoiding involution. They are
enumerated by $a_{n-1} = {n-1 \choose {[(n-1)/2]}}$.
\\
By $(3)$ of Remark~\ref{rmq-ssd} we deduce when $n$ and $s$ are even
that for $\pi=\pi'\pi''$ with $|\pi'|=|\pi''|$ and
$\pi''\in\SS_{n/2}(1\mn3\mn2)$ having $s/2$ right-to-left maxima.
This mapping is trivially bijective. Moreover, the number of
$1\mn3\mn2$-avoiding permutations according to the length and to
the number of right-to-left maxima is given by the ballot numbers
(see the sorting algorithm with one stack \cite{Kn} described in
the end of Subsection~\ref{intro-tools}).

So, the formula is established from a combinatorial way for the case
of an even number of right-to-left maxima and for the special case
of one right-to-left maximum.
Two simple bijections given in Lemma~\ref{lemm-ssd1} and Lemma~\ref{lemm-ssd2}
allow to compute for the case of an odd (and greater than $1$) number
of right-to-left maxima,
but we now state its formula from a combinatorial way.

First, we prove that all the arguments established for
$1\mn3\mn2$-avoiding involutions hold for Dyck word prefixes that
is the number involutions in $\II_n(1\mn3\mn2)$ having $2l$ or
$2l+1$ right-to-left maxima is equal to the number of Dyck word
prefixes of length $n$ ended by $x \ax^l$ which are also Dyck
words or not respectively. Of course, a Dyck word prefix of odd
length cannot be also a Dyck word. Let $w = w' x \overline{x}^l$
be a Dyck word prefix of length $n$. Trivially, if $l=0$, $w=w'x$
is in bijection with $w'$ a Dyck word prefix of length $n-1$. It
is well known that all Dyck words of length $n=2k$ ended by $x
\ax^l$ are enumerated by the ballot number ${2k-l-1 \choose k-1} -
{2k-l-1 \choose k}$. Moreover, if $w$ is not also a Dyck word then
$w$ is in bijection with a Dyck word prefix $w\ax$ of length $n+1$
such that if $n=2k$ then $w\ax$ is not also a Dyck word whereas if
$n=2k+1$ then $w\ax$ is also a Dyck word if and only if
$|w|_x=|w|_{\ax}+1$.

We now complete the proof by establishing that the number of Dyck
word prefixes $w = w' x \ax^l$ of length $2k+1$ (which are not
Dyck words) is ${2k-l \choose k,k-l}$. Let $w' = u x v$ such that
$v x \ax^l$ is a Dyck word (so of even length, thus like $u$). Let
$u'$ be the bilateral word in bijection with $u$ by $\Xi$. Thus,
$u' x v$ is any word of $\{x,\ax\}^*$ of length $2k-l$ with $k$
letters $x$ (and $k-l$ letters $\ax$). Conversely, let $w''$ be a
word of $\{x,\ax\}^*$ with $|w''|=2k-l$ and $|w''|_x=k$ (and
$|w''|_{\ax}=k-l$). There exists a nonnegative integer $m$ such
that $w'' = w_0 \overline{x} w_1 \ax \ldots \ax w_{m} x w_{m+1} x
\ldots x w_{2m+l}$ where $w_i$ is a Dyck word for all $0 \leq i
\leq 2m+l$. Let $u' = w_0 \overline{x} w_1 \ax \ldots \ax w_{m} x
w_{m+1} x \ldots x w_{2m}$ and $v = w_{2m+1} x w_{2m+2} x \ldots x
w_{2m+l}$. So, by applying bijection $\Xi$ on $u'$, we have $u =
w_0 x w_1 x \ldots x w_{m} x w_{m+1} x \ldots x w_{2m}$. Thus, we
obtain $w = u x v x \ax^l$.\qed

\subsection{Counting an occurrences of $12\mn3\mn\dots\mn k$}

Now, let us define
$$\begin{array}{l}
D_I(x_1,x_2,\dots)=\sum\limits_{\pi\in\mathcal{I}(1\mn3\mn2)}
x_1^{\sharp1(\pi)}\prod\limits_{k\geq 2}
x_k^{\sharp12\mn3\mn\dots\mn k(\pi)},\\
D_S(x_1,x_2,\dots)=\sum\limits_{\pi\in\mathcal{S}(1\mn3\mn2)}x_1^{\sharp1(\pi)}\prod\limits_{k\geq
2} x_k^{\sharp12\mn3\mn\dots\mn k(\pi)}. \end{array}$$

\begin{theorem}\label{thd}
The generating function $D_I(x_1,x_2,\dots)$ is given by
$$D_I(x_1,x_2,\dots)=\dfrac{1+x_1-x_1x_2+x_1x_2D_I(x_1,x_2x_3,x_3x_4,\dots)}
        {1-x_1^2+x_1^2x_2^2-x_1^2x_2^2D_S(x_1^2,x_2^2x_3^2,x_3^2x_4^2,\dots)},$$
where {\rm(}see \cite[Theorem~1]{MS}{\rm)}
$$D_S(x_1,x_2,x_3,x_4,\dots)=\dfrac{1}{1-x_1+x_1x_2-x_1x_2D_S(x_1,x_2x_3,x_3x_4,\dots)}.$$
\end{theorem}
\begin{proof}
In \cite[Theorem~1]{MS} proved
    $$D_S(x_1,x_2,x_3,x_4,\dots)=\frac{1}{1-x_1+x_1x_2-x_1x_2D_S(x_1,x_2x_3,x_3x_4,\dots)}.$$
On the other hand, by Proposition \ref{prom} we have exactly two
block decompositions for an arbitrary $1\mn3\mn2$-avoiding involution.
The contribution of the first case is, if $n=1$ we get $x_1$,
otherwise we have $x_1x_2(D_I(x_1,x_2x_3,x_3x_4,\dots)-1)$. The
contribution of the second case, if $\beta$ is empty (see
Proposition~\ref{prom}) we get $x_1^2D_I(x_1,x_2,\dots)$,
otherwise we have (see Lemma~\ref{genl})
$x_1^2x_2^2(D_S(x_1^2,x_2^2x_3^2,x_3^2x_4^2,\dots)-1)D_I(x_1,x_2,\dots)$.
Therefore
$$\begin{array}{l}
D_I(x_1,x_2,\dots)=1+x_1-x_1x_2+x_1x_2D_I(x_1,x_2x_3,x_3x_4,\dots)+x_1^2D_I(x_1,x_2,\dots)+\\
\qquad\qquad\qquad\qquad\qquad\qquad\qquad\qquad\qquad\qquad+x_1^2x_2^2D_I(x_1,x_2,\dots)(D_S(x_1^2,x_2^2x_3^2,\dots)-1).
\end{array}$$
The rest is easy to check.
\end{proof}

\begin{example}
By \cite[Proposition~3.1]{MS} we get
$$f(x,y):=D_S(x^2,y^2,1,1,\dots)=\frac{-1+x^2-\sqrt{1-2x^2+x^4+4x^2y^2}}{2x^2y^2}.$$
So by Theorem~\ref{thd} we have
$$\sum_{\pi\in\mathcal{I}(1\mn3\mn2)}x^{|\pi|}y^{\sharp12(\pi)}=
\frac{1+x-xy}{1-x^2-xy+x^2y^2-x^2y^2f(x,y)}.$$ For example, the
number of involutions in $\II_n(1\mn3\mn2)$ containing exactly
once a generalized pattern $12$ is given by
$\frac{1}{4}(2n-1+(-1)^n)$, and the number of involutions in
$\II_n(1\mn3\mn2)$ containing exactly twice a generalized pattern
$12$ is given by $\frac{1}{8}(2n(n-2)+1-(-1)^n)$.
\end{example}

Of course this last result for
$\sum_{\pi\in\mathcal{I}(1\mn3\mn2)}x^{|\pi|}y^{\sharp12(\pi)}$ is
equivalent to Theorem~\ref{rises} enumerating the number of
involutions in $\II_n(1\mn3\mn2)$ having $r$ rises (or
equivalently having $n-r$ left-to-right minima by
Remark~\ref{lrm=rises}) with $|\pi|=n$ and with $\sharp12(\pi)=r$.

\end{document}